\documentclass{mcom-l}

\newcommand{\R}{\mathbb{R}}
 
\newcommand{\Z}{\mathbb{Z}}
\newcommand{\N}{\mathbb{N}}

\newcommand{\PP}{\mathbb{P}}

\DeclareMathOperator{\ord}{ord}
\DeclareMathOperator{\lcm}{lcm}
\DeclareMathOperator{\integ}{int}

\usepackage{amssymb}
\usepackage{algpseudocode}
\algtext*{EndWhile}
\algtext*{EndFor}
\algtext*{EndIf}

\copyrightinfo{}{}

\newtheorem{theorem}{Theorem}[section]
\newtheorem{lemma}[theorem]{Lemma}
\newtheorem{corollary}[theorem]{Corollary}

\theoremstyle{definition}
\newtheorem{definition}[theorem]{Definition}

\newtheorem{problem}[theorem]{Problem}
\newtheorem{algorithm}[theorem]{Algorithm}

\theoremstyle{remark}
\newtheorem{remark}[theorem]{Remark}

\numberwithin{equation}{section}

\begin{document}

\title[BSGS Integer factorization]{A babystep-giantstep method for faster deterministic integer factorization}


\author[M. Hittmeir]{Markus Hittmeir}
\address{Hellbrunnerstra{\ss}e 34, A-5020 Salzburg}
\curraddr{}
\email{markus.hittmeir@sbg.ac.at}
\thanks{The author is supported by the Austrian Science Fund (FWF): Project F5504-N26,
which is a part of the Special Research Program "Quasi-Monte Carlo Methods: Theory
and Applications".}


\date{}

\dedicatory{}

\begin{abstract}
In 1977, Strassen presented a deterministic and rigorous algorithm for solving the problem of computing the prime factorization of natural numbers $N$. His method is based on fast polynomial arithmetic techniques and runs in time $\widetilde{O}(N^{1/4})$, which has been state of the art for the last forty years. In this paper, we will combine Strassen's approach with a babystep-giantstep method to improve the currently best known bound by a superpolynomial factor. The runtime complexity of our algorithm is of the form
\[
\widetilde{O}\left(N^{1/4}\exp(-C\log N/\log\log N)\right).
\]
\end{abstract}

\maketitle

\section{Introduction}
Let $N$ be any natural number. We consider the problem of computing the prime factorization of $N$. In practice, a large variety of probabilistic and heuristic methods is used for this  task. We refer the reader to the survey \cite{Len} and to the monographs \cite{Rie} and \cite{Wag}. In this paper, we discuss a more theoretical aspect of the integer factorization problem and focus on \emph{deterministic} algorithms and their rigorous analysis. We will describe runtime complexities by using the bit-complexity model of the multitape Turing machine. 

In \cite{Str}, Strassen used fast polynomial multiplication and multipoint evaluation techniques to establish a deterministic and rigorous factorization algorithm running in $O(N^{1/4+o(1)})$ bit operations. Since the publication of Strassen's approach in 1977, there have been a few refinements of the logarithmic factors in the bound for the runtime. Currently, the best known bit-complexity bound is given by
\[
O\Big{(}{\textsf{M}}_{\integ}\Big{(}\frac{N^{1/4}\log N}{\sqrt{\log\log N}}\Big{)}\Big{)}
\]
and has been proven by Costa and Harvey in \cite{CosHar}. ${\textsf{M}}_{\integ}(k)$ denotes the cost of multiplying two $\lceil k\rceil$-bit integers and is in $O(k\log k\cdot 8^{\log^*k})$, where $\log^*k$ is the iterated logarithm (see \cite{HarHoeLec}). Costa and Harvey improved Strassen's method and refined the preceding result of Bostan, Gaudry and Schost \cite{BosGauSch} by the factor of $\sqrt{\log\log N}$. The goal of our paper is to establish a new approach for deterministic integer factorization. The main result improves the currently best known runtime complexity bound by a factor growing faster than \emph{any fixed power} of $\log N$.

To describe our approach, we consider the case $N=pq$, where $p$ and $q$ are unknown prime factors and $p<q$. We will employ a refined babystep-giantstep method to solve the discrete logarithm problem $a^X\equiv a^{N+1}\mod N$ for a certain $a\in\Z$ coprime to $N$. The purpose of this procedure is to determine $S:=p+q$, which is a solution to this problem due to Euler's theorem. Knowing $S$ allows us to factor $N$ immediately. In our main algorithm, we first compute a lower bound for $p$ with Strassen's approach and then use the resulting upper bound for $S$ in the babystep-giantstep method.

In addition to the preceding observations, let $m\in\N$. We will consider so called modular hyperbolas $\mathcal{H}_{N,m}$, which are defined as the sets of solutions $(x,y)$ to the congruence equation $N\equiv xy \mod m$. Clearly, the corresponding set $\mathcal{L}_{N,m}$ consisting of the elements $x+y\mod m$ for $(x,y)\in\mathcal{H}_{N,m}$ contains the residue of $S$ modulo $m$.
If $r$ is prime, than the cardinality of $\mathcal{L}_{N,r}$ is about half of the possible residue classes modulo $r$. Considering all primes up to a suitable bound $B$, we deduce significant information about $S$. 

We discuss a brief example. Let $N:=7909787=3823\cdot 2069$. For $r=5$, we have $N\equiv 2 \mod 5$ and $(1,2),(2,1),(3,4),(4,3)\in\mathcal{H}_{N,5}$. Consequently, $2,3\in\mathcal{L}_{N,5}$. Considering all primes up to $11$ and putting $m=2\cdot 3\cdot 5\cdot 7\cdot 11$, we observe that $\mathcal{L}_{N,m}$ contains only $40$ elements. As a consequence, $S\mod m$ is restricted to $40/2310=1.7\%$ of the residue classes modulo $m$. We may use this fact in our refined babystep-giantstep method for determining $S$ itself. This application yields the main contribution to our improvement.

In our analysis of the logarithmic factors in the runtime complexity of our algorithm, we will make use of the following notion.
\definition{We call a function $\textsf{e}_{\integ}:\R\rightarrow\R$ \emph{logarithmic multiplication time}, if
\begin{enumerate}
\item{$\lceil k\rceil$-bit integers can be multiplied in at most $O(k\cdot\textsf{e}_{\integ}(k))$ bit operations,}
\item{$\textsf{e}_{\integ}(k^{1/c})\in \Theta(\textsf{e}_{\integ}(k))$ for every constant $c\in\N$,}
\item{$\textsf{e}_{\integ}(k)\in \Omega(\log k)$.} 
\end{enumerate}\label{d1}}

Considering the mentioned bound for ${\textsf{M}}_{\integ}(k)$, one observes that $\log k\cdot 8^{\log^*k}$ is a logarithmic multiplication time. We are now ready to formulate our main theorem.

{\theorem{Let $N$ be any natural number. Then the prime factorization of $N$ can be computed in
\[
O\left(N^{\frac{1}{4}}\exp\left(\frac{-C_1\log N}{\log \left(C_2\log N\right)}\right)\cdot\sqrt{\emph{\textsf{e}}_{\integ}(N)\log^3 N}\right)
\]
bit operations, where $C_1=(\log 2)/16$ and $C_2=0.2103$.}\label{A}}
\vspace{12pt}

We give a brief overview on the contents of this paper:

Section 2: We discuss auxiliary results required in later sections of the paper.

Section 3: We adapt recent results concerning Strassen's approach for usage in our methods, especially for the main algorithm in Section 8.

Section 4: Let $m_1,\ldots ,m_d\in\N$ be pairwise coprime, let $M$ be their product and $\mathcal{A}_i\subseteq \Z_{m_i}$ for $i=1,\ldots,d$. By the Chinese Remainder Theorem applied to the system of congruences modulo the $m_i$, we know that there is an unique set $\mathcal{A}\subseteq \Z_{M}$ determined by the $\mathcal{A}_i$. In this section, we present an \emph{asymptotically optimal} algorithm for computing the elements of $\mathcal{A}$. This procedure is used later on to efficiently compute the sets $\mathcal{L}_{N,m}$ mentioned above.

Section 5: Here, the key ingredients for the proof of Theorem \ref{A} are developed, including the study of the modular hyperbolas and the bounds of the cardinality of the sets $\mathcal{L}_{N,m}$.

Section 6: Two algorithms are discussed in this section. The first one, Theorem \ref{t2}, is due to Sutherland \cite{Sut} and improves Shanks' babystep-giantstep technique for the problem of finding the order of elements modulo $N$ by the factor of $\sqrt{\log\log N}$. The second, Algorithm \ref{a1}, makes use of Theorem \ref{t2}. It is a sufficiently fast solution to the following problem: Given a natural number $N$ and some $\delta\in\R^+$, either find a nontrivial factor of $N$, prove that $N$ is prime, or find an element in its multiplicative group whose order is larger than $\delta$.

Section 7: Theorem \ref{t3}, one of the key results for the proof of Theorem \ref{A}, is shown. This statement concerns the problem of finding a set of solutions to $a^X\equiv a^{N+1}\mod N$ containing $S$.

Section 8: The results from previous sections are combined to prove our main result. The runtime complexity is analyzed and we also show that Theorem \ref{A} is a corollary to the slightly stronger statement in Theorem \ref{t9}.

\section{Elementary and analytical results}
This section contains auxiliary results. We will discuss elementary lemmas, Shank's babystep-giantstep method for solving the discrete logarithm problem and some analytical results about the prime-counting function and Chebyshev's theta function. All of these are either well-known or easy to prove. 

Throughout the paper, $\PP$ denotes the set of prime numbers. We call a natural number \emph{semiprime} if and only if it is the product of two distinct primes. Furthermore, let $\N_0:=\N\cup\{0\}$ and $\Z_n:=\Z/n\Z$ for $n\in\N$. 

{\lemma{Let $u,v\in\N$ such that $u\leq v$ and assume that we know $N:=uv$ and $S:=u+v$. Then we can compute $u$ and $v$ in $O(\emph{\textsf{M}}_{\integ}(\log N))$  bit operations.}\label{l1}}

\begin{proof}
Set $D:=v-u$ and note that $S^2-D^2=4N$. Hence, $D$ is the positive solution of $\sqrt{S^2-4N}$, and we derive $v=(S+D)/2$ and $u=N/v$. Since squaring $S$ is most expensive, the number of required bit operations for these computations is bounded by $O(\textsf{M}_{\integ}(\log N))$ .
\end{proof}

{\lemma{Let $N$ be semiprime with prime factors $p$ and $q$. If $a\in\Z$ is coprime to $N$, then $a^{p+q}\equiv a^{N+1}\mod N$.}\label{l3}}
\begin{proof}
This is an immediate consequence of Euler's Theorem and the fact that $\phi(N)=(p-1)(q-1)=N+1-(p+q)$.
\end{proof}

{\lemma{Let $N\in\N$ and $a\in\Z$ be coprime. Furthermore, define $m:=\ord_N(a)$. Then
\[
\ord_p(a)=m \text{\hspace{12pt} for all primes $p$ dividing $N$}
\]
if and only if
\[
\gcd(N,a^{m/r}-1)=1 \text{\hspace{12pt} for all primes $r$ dividing $m$}.
\]\label{l2}}
\vspace{-1em}
\begin{proof}
We show that the first statement implies the second. Let $r$ be an arbitrary prime dividing $m$. Then $a^{m/r}\not\equiv 1 \mod p$ for every prime $p$ with $\ord_p(a)=m$. But since this holds for every prime factor of $N$, we derive $\gcd(N,a^{m/r}-1)=1$. The prime divisor $r$ of $m$ was arbitrary, hence, the statement follows.

For proving the reverse statement, let $p$ be an arbitrary prime dividing $N$. We want to show that $\ord_p(a)=m$. Clearly, we have $a^{m}\equiv 1 \mod p$. Therefore, the order of $a$ modulo $p$ is a divisor of $m$. Assume to the contrary that $\ord_p(a)$ is a proper divisor $d$ of $m$. Then $a^{d}\equiv 1 \mod p$. This implies that some prime divisor $r$ of $m$ exists such that $a^{m/r}\equiv 1 \mod p$, which contradicts our assumption.
\end{proof}}

Sections 6 and 7 focus on finding solutions to certain cases of the following computational problem. We are interested in unconditional and deterministic runtime complexity bounds. The following remark discusses the state of the art, the so called babystep-giantstep method by Shanks. In the proofs of the results of the Sections 6 and 7, we will consider several variants of this idea.

{\problem{(Constrained Discrete Logarithm Problem, CDLP)

Let $N\in\N$, $1\leq E<T\leq N$ and $a,b\in\Z_N^*$. Find the smallest $x\in [E,T]$ such that 
$
a^x\equiv b \mod N,
$
or prove that no such $x$ exists.\label{p1}}}

\remark{(Shanks' algorithm)

We recall Shanks' babystep-giantstep method, which is a deterministic algorithm to solve this problem performing 
\[
O(\sqrt{T-E}\cdot \log^2 N+\log E\cdot \textsf{M}_{\integ}(\log N))
\]
bit operations. Let $m:=\lceil\sqrt{T-E}\rceil$ and write $x=mu+v+E$ for unknown $u,v\in\{0,\ldots,m-1\}$. We first compute $a^{i}a^{E}\mod N$ for $i=0,\ldots,m-1$, the so called $\textit{babysteps}$, which requires $m+O(\log E)$ multiplications. We proceed by computing $a^{-mj}b\mod N$ for $j=0,\ldots,m-1$, the so called $\textit{giantsteps}$. If we first precompute $a^{-m}\mod N$ and then use the recursive formula
\[
a^{-m(j+1)}b\equiv a^{-m}a^{-mj}b \mod N,
\]
this needs $O(m)$ multiplications. We now want to check if one of the giantsteps is equal to one of the babysteps. In the multitape Turing machine model, this can be done by using merge sort. We denote every number as string of at most $\lceil \log N \rceil$ bits and sort both the list of babysteps and the list of giantsteps by performing at most $O(m\log m)$ comparisons, yielding a bit-complexity bounded by $O(m\log^2 N)$. After applying merge sort, it is easy to see that one is able to find all matches in the lists by $O(m\log N)$ bitwise comparisons. For more detailed information on merge sort, we refer the reader to \cite[Chap. 2.2]{SedWay}. For the claimed complexity, see Proposition F in \cite[Page 272]{SedWay}. 

To summarize, we need $O(m+\log E)$ group operations for computing the baby- and giantsteps and at most $O(m \log^2 N)$ bit operations for the sorting and searching process. Since one multiplication in $\Z_N^*$ is bounded by $O(\textsf{M}_{\integ}(\log N))$ bit operations and $\textsf{M}_{\integ}(\log N)< \log^2 N$, we deduce the stated running time.
\label{r1}}

\definition{Let $x\in\R^+$. We define the prime counting function
\[
\pi(x):=\#\{r\leq x: r\in\PP\}
\]
and Chebyshev's theta function
\[
\vartheta(x):=\sum_{r\in\PP:\text{ } r\leq x}\log r.
\]
}
The Prime Number Theorem states that $\pi(x)\sim x/\log x$, which also implies that $\vartheta(x)-x\in o(x)$. In the following theorem, we consider some explicit estimates for both $\pi(x)$ and $\vartheta(x)$.

{\theorem{Let $x\in\R^+$. We have
\begin{enumerate}
\item{$x/\log x<\pi(x)$ for $x\geq 17$.}
\item{$\vartheta(x)-x<x/(4\log x)$ for $x>1$.}
\end{enumerate}
}\label{t8}}
\begin{proof}
The first statement is Corollary 1 on page 69 of \cite{RosSch}. The second statement is an immediate consequence of Inequality 5.6a* of Theorem 7* in \cite{Sch}.
\end{proof}

In Theorem \ref{t3}, one of the key results of this paper, we will employ the following inequality.

{\lemma{Let $x\in\R^+$, $x>1$. Then
\[
\vartheta(x)-x<\frac{\log 2}{2}\pi(x)-\frac{1}{2}\log\log x.
\]}\label{l7}}

\begin{proof}
It is easy to show that, for every $x\geq 20$, we have 
\[
\frac{1}{2}\log\log x <\frac{x}{12\log x}.
\]
Combining this with the first statement of Theorem \ref{t8}, we deduce that
\[
\frac{x}{4\log x}<\left(\frac{\log 2}{2}-\frac{1}{12}\right)\frac{x}{\log x}<\frac{\log 2}{2}\pi(x)-\frac{1}{2}\log\log x
\]
for every $x\geq 20$. From the second statement of Theorem \ref{t8}, it follows that the claim is true for $x\geq 20$. Since one can easily verify the claim for $x<20$, this finishes the proof.
\end{proof}

{\theorem{
For primes $r$, we have 
\[
\prod_{r\leq B}\frac{r-1}{r}\in\Theta\Big{(}\frac{1}{\log B}\Big{)}.
\]}}

\begin{proof}
This result is called Mertens' Theorem. A proof can be found in \cite[Thm. 5.13]{Sho}.
\end{proof}

In Section \ref{S2}, we will need the following statement. It is equivalent to Mertens' Theorem.

{\corollary{
For primes $r$, we have
\[
\prod_{r\leq B}\frac{r+1}{r}\in\Theta(\log B).
\]
}\label{t5}}

\begin{proof}
The claim follows from Mertens' Theorem and the fact that
\[
\left(\prod_{r\leq B}\frac{r+1}{r}\right)\left(\prod_{r\leq B}\frac{r-1}{r}\right)=\prod_{r\leq B} \left(1-\frac{1}{r^2}\right)\in\Theta(1).
\]
\end{proof}

\section{Strassen's approach}
Let us briefly recall Strassen's idea for factoring natural numbers. Let $N\in\N$ and set $d:=\lceil N^{1/4}\rceil$. We want to compute subproducts of the product $\lceil N^{1/2}\rceil !$ to find a factor of $N$. For this task, the polynomial 
\[
g=(X+1)(X+2)\cdots (X+d)
\]
is computed modulo $N$ and evaluated at the points $0,d,2d,\ldots,(d-1)d$ by using fast polynomial arithmetic techniques. Next, $g_i:=\gcd(g(id),N)$ is computed for $i=0,\ldots,d-1$. If $N$ is not a prime, one of these GCDs is nontrivial. If $g_{i_0}=N$, we obtain a nontrivial factor of $N$ by computing $\gcd(i_0d+j,N)$ for $j=1,\ldots,d$.

As we have already mentioned in the introduction, we will make use of recent results that are closely related to Strassen's approach. First, we will employ a theorem of \cite{Hit} in Algorithm \ref{a1}. Second, we will adapt a proposition of \cite{CosHar} in Algorithm \ref{a2}, which is our main algorithm for factoring semiprime numbers.

{\theorem{Let $N$ be a natural number  and $s,m\in\N_0$ such that $s,m<N$ and $s\equiv p\mod m$ for every prime divisor $p$ of $N$. Knowing $s$ and $m$, one can compute the prime factorization of $N$ in 
\[
O\Big{(}\emph{\textsf{M}}_{\integ}\Big{(}\frac{N^{1/4}}{\sqrt{m}}\log N\Big{)}\Big{)}
\]
bit operations.}\label{t4}}

\begin{proof}
This is Theorem 2.8 in \cite{Hit}. In the proof, we use another polynomial instead of $g$, taking the additional information about the prime factors of $N$ into account. Essentially, we are getting rid of all the factors of $\lceil N^{1/2}\rceil !$ which are not of the form $mj+s$, $j\in\N$.
\end{proof}

Now let $B\in\R$ and put $Q=\prod_{r\leq B} r$ for prime numbers $r$. Furthermore, set $\rho=\phi(Q)=\prod_{r\leq B}(r-1)$. For the proof of Theorem \ref{A}, we adapt the following result of Costa and Harvey.

{\theorem{Let $k\geq 0$ and $b=4^k\rho Q$. Assume that $b<N$ and $\gcd(N,Q)=1$. We can find a prime divisor $q$ of $N$ such that $q\leq b$, or prove that no such divisor exists, in
\[
O\Big{(}\emph{\textsf{M}}_{\integ}(2^k\rho\log N)+(Q^2+\log(2^k\rho))\emph{\textsf{M}}_{\integ}(\log N)\log\log N\Big{)}
\]
bit operations.}}

\begin{proof}
This is Proposition 9 in \cite{CosHar}. The proof is based on the observation that if $\gcd(N,Q)=1$, all multiples of the primes smaller or equal to $B$ cannot be factors of $N$. Again, we replace $g$ by another polynomial to get rid of these numbers, which leads to the savings in the runtime complexity.
\end{proof}

{\corollary{Let $\delta\in\R$ such that $N^{1/3}\leq \delta\leq N/4$. We can find a prime divisor $q$ of $N$ such that $q\leq \delta$, or prove that no such divisor exists, in
\[
O\Big{(}\emph{\textsf{M}}_{\integ}\Big{(}\frac{\sqrt{\delta}\log N}{\sqrt{\log\log N}}\Big{)}\Big{)}
\]
bit operations.\label{c3}}}

\begin{proof}
Let $B=\frac{1}{19}\log N$ and $Q$, $\rho$ as above. For Chebyshev's theta function $\vartheta(x)$, the Prime Number Theorem yields that $\vartheta(x)-x\in o(x)$ and therefore
\[
Q=\prod_{r\leq B}r=e^{\vartheta(B)}\in O(e^{B+o(B)})\subseteq O(N^{1/18}).
\]
Note that $\rho Q<Q^2\in O(N^{1/9})$ and let $k\geq 0$ such that $4^{k-1}\rho Q<\delta \leq 4^{k}\rho Q$. After trial division up to $B$, we either find a prime divisor or we have proven that $\gcd(N,Q)=1$. Assume the latter. We now apply the previous theorem with the value $b:=4^k\rho Q <4\delta \leq N$. Note that $b\in O(\delta)$ and therefore we have
$2^k\rho=\sqrt{b\rho/Q}\in O\Big{(}\sqrt{\delta\rho/Q}\Big{)}$. Using Mertens' Theorem, we derive that
\[
\frac{\rho}{Q}=\prod_{r\leq B}\frac{r-1}{r}\in O\left(\frac{1}{\log\log N}\right).
\]
Hence, the first summand in the running time of the theorem is of the claimed form. The second summand is in $O((Q^2)^{1+o(1)})\subseteq O(N^{1/9+o(1)})$, but since we have $\sqrt{\delta}\geq N^{1/6}$, this is asymptotically negligible.
\end{proof}

\section{Optimal Chinese Remaindering\label{S3}}
Let $m_1,\ldots,m_d\in\N$ be pairwise coprime and $\mathcal{A}_i\subseteq \Z_{m_i}$ for $i=1,\ldots,d$. Furthermore, define $M:=\prod_{i=1}^d m_i$. By the Chinese Remainder Theorem, we know that there is a unique set $\mathcal{A}\subseteq\Z_M$ determined by the sets $\mathcal{A}_i$. We may define this set as
\[
\mathcal{A}:=\left\{x\in\Z_M\text{ }\Bigg{|}\text{ } \exists (\alpha_1,\ldots,\alpha_d)\in\prod_{i=1}^d \mathcal{A}_i:x\equiv \alpha_i\mod m_i,\text{ } 1\leq i \leq d\right\}.
\]
In this section, we will discuss an asymptotically optimal algorithm for determining the set $\mathcal{A}$. Ignoring implied constants, the method takes the same time to compute all its elements as it would take to simply write them down. For $i=1,\ldots,d$, we set $\kappa_i:=|\mathcal{A}_i|$ and $\mathcal{A}_i=\{a_{i,1},\ldots,a_{i,\kappa_i}\}$. Furthermore, let $M_i:=Mc_i/m_i$, where $c_i:=(M/m_i)^{-1} \mod m_i$. With this notation, consider the following procedure.

\begin{algorithm}
\emph{Input:} Pairwise coprime $m_i\in\N$, $\mathcal{A}_i\subseteq \Z_{m_i}$ and $\kappa_i$ for $i=1,...,d$, where $\kappa_1\geq \kappa_2\geq\ldots\geq \kappa_d$.

\emph{Output:} The set $\mathcal{A}\subseteq\Z_M$ determined by the $\mathcal{A}_i$ due to the CRT.

\begin{algorithmic}[1]
\State Compute $M$ and $|\mathcal{A}|=\prod_{i=1}^d \kappa_i$. 
\For {$i=1,\ldots,d$} \Comment{Precomputation  loop}
	\State Compute $M_i\mod M$ and set $r_i=0$. 
	\State $\Delta_{i,0}\gets (a_{i,1}-a_{i,\kappa_i})M_i\mod M$
	\For{$l_i=1,\ldots,\kappa_i-1$}
		\State $\Delta_{i,l_i}\gets (a_{i,l_i+1}-a_{i,l_i})M_i \mod M$ 
	\EndFor
\EndFor
\State $x_1\gets a_{1,1}M_1+...+a_{d,1}M_d \mod M$
\For {$\nu=1,\ldots,|\mathcal{A}|-1$} \Comment{Main loop}
		\State Compute $r_1\gets r_1+1\mod{\kappa_1}$ and set $\mu=1$.
		\While {$r_\mu=0$}
			\State $\mu\gets \mu+1$
			\State $r_\mu\gets r_\mu+1\mod{\kappa_\mu}$
		\EndWhile
		\State $x_{\nu+1}\gets x_{\nu}+\sum_{i=1}^\mu\Delta_{i,r_i} \mod M$
\EndFor
\end{algorithmic}
\label{a3}
\end{algorithm}

{\theorem{Algorithm \ref{a3} is correct. Assuming $|\mathcal{A}|\geq \sum_{i=1}^d \kappa_i\cdot \log M$, its runtime complexity is bounded by $O(|\mathcal{A}|\cdot \log M)$ bit operations.}}

\begin{proof}
\textit{Correctness:} We follow Gauss' approach to find the elements of $\mathcal{A}$. With our notation, we have to compute the $|\mathcal{A}|$ elements
\begin{equation}
a_{1,l_1}M_1+a_{2,l_2}M_2+\cdots+a_{d,l_d}M_{d} \mod M,
\end{equation}
where $l_i\in\{1,\ldots,\kappa_i\}$. Steps 1-7 are precomputational. To verify correctness, we discuss the idea behind the main loop from Step 8 to 13. We start with the value $x_1=a_{1,1}M_1+...+a_{d,1}M_d \mod M$, which has already been computed in Step $7$. We now change the value $a_{1,1}$ to $a_{1,2}$ in this representation. It is easy to see that we may do this by simply computing $x_2=x_1+\Delta_{1,1}\mod M$. Next, we change $a_{1,2}$ to $a_{1,3}$ by computing $x_3=x_2+\Delta_{1,2}\mod M$. We proceed in this way until we reach $x_{\kappa_1}=x_{\kappa_1-1}+\Delta_{1,\kappa_1-1}\mod M$. In the next run of the loop, we have $r_1=0$. Consequently, we do not leave the small loop for the values of $\mu$, and $a_{2,1}$ changes to $a_{2,2}$. After incrementing $r_2$, we get
\[
x_{\kappa_1+1}=x_{\kappa_1}+\Delta_{1,0}+\Delta_{2,1} \mod M.
\]
Note that adding $\Delta_{1,0}$ to $x_{\kappa_1}$ brings back the value $a_{1,\kappa_1}$ to $a_{1,1}$ in the representation of $x_{\kappa_1}$ as a sum of the form (4.1). We now proceed by adding the values $\Delta_{1,l_1}$ again and compute $x_{\kappa_1+1},\ldots,x_{2\kappa_1}$, where we have to change $a_{2,2}$ to $a_{2,3}$. As soon as we have reached $a_{2,\kappa_2}$, the value $a_{3,1}$ has to be changed to $a_{3,2}$. Imagine this procedure as a clock, where the values of $a_{1,l_1}$ represent a small hand and the values of $a_{2,l_2}$ represent a bigger hand which is moved every time the small hand has made a complete round. The values of $a_{3,l_3}$ represent an even bigger hand, and so on. After $|\mathcal{A}|-1$ runs of the main loop, we have computed all the elements of $\mathcal{A}$.

\textit{Running time:}
The cost for Step 1 is bounded by $O(d\cdot \textsf{M}_{\integ}(\log M))$ bit operations. Since our assumption implies 
$
|\mathcal{A}|\log M\geq \sum_{i=1}^d \kappa_i\cdot \log^2 M\geq d\cdot \textsf{M}_{\integ}(\log M),
$
this is negligible. In each run of the loop from Step 2 to 6, we have to perform $O(\kappa_i)$ multiplications modulo $M$. Consequently, the total runtime of this loop is in $O(\sum_{i=1}^d \kappa_i \cdot \textsf{M}_{\integ}(\log M))$ bit operations, which is also negligible. 

We are left with discussing the complexity of the main loop. Considering the precomputation loop, the values $\Delta_{i,j}$ are stored on tape in the following order: $\Delta_{1,0},\Delta_{1,1},...,\Delta_{1,\kappa_1-1}$, followed by $\Delta_{2,0},...,\Delta_{2,\kappa_2-1}$ and so on. There are $|\mathcal{A}|-1$ runs of the main loop. In each of these runs, we have to increment $r_1$ modulo $\kappa_1$ and add $\Delta_{1,r_1}$ to $x_\nu$. For doing so, we may simply walk along the list of the $\Delta_{1,j}$. Since addition of numbers smaller than $M$ is in $O(\log M)$ bit operations, the cost for these calculations is $O(|\mathcal{A}|\log M)$. Furthermore, in every $\kappa_1$-th run, we additionally have to increment the value of $r_2$ modulo $\kappa_2$ and add $\Delta_{2,r_2}$ to $x_\nu$. In these cases, we also have to account for the cost of seeking to the correct value in the list $\Delta_{2,0},...,\Delta_{2,\kappa_2-1}$ and going back to $\Delta_{1,0}$, which is stored at the beginning of the tape. One observes that this can be done by $O((\kappa_1+2\kappa_2)\log M)$ bit operations. Since there are at most $|\mathcal{A}|/\kappa_1$ such runs, we may bound the cost by
\[
O\left(\frac{\kappa_1+2\kappa_2}{\kappa_1}\cdot|\mathcal{A}|\log M\right).
\]
Summing up, the overall bit-complexity bound for the main loop is given by
\[
O\left(\left(1+\sum_{i=1}^{s-1}\frac{\kappa_1+2(\kappa_2+\cdots+\kappa_{i+1})}{\kappa_1\cdots \kappa_i}\right)\cdot|\mathcal{A}|\log M \right).
\]
Here, $s$ denotes the largest number such that $\kappa_s>1$. Note at this point that we have $\kappa_1\geq \kappa_2\geq...\geq \kappa_d$. Hence, we derive that
\[
\sum_{i=1}^{s-1}\frac{\kappa_1+2(\kappa_2+\cdots+\kappa_{i+1})}{\kappa_1\cdots \kappa_i}\leq \sum_{i=1}^{s-1}\frac{(1+2i)\kappa_1}{\kappa_1\cdots \kappa_i}\leq \sum_{i=1}^{\infty}\frac{1+2i}{2^{i-1}}
\] 
is constant, which proves the claim.
\end{proof}

{\corollary{There exists an algorithm with the following properties:
\begin{itemize}
\item{Input: Pairwise coprime $m_i\in\N$, $\mathcal{A}_i\subseteq \Z_{m_i}$ and $\kappa_i$ for $i=1,...,d$, where $\kappa_1\geq \kappa_2\geq\ldots\geq \kappa_d$. Furthermore, a natural number $N$ and $a\in\Z_N^*$.}
\item{Output: The elements $a^{x_\nu}\mod N$ for $x_\nu \in \mathcal{A}$.}
\end{itemize}
Assuming $|\mathcal{A}|\geq \sum_{i=1}^d \kappa_i\cdot \log M$, the runtime complexity of the algorithm is bounded by $O(|\mathcal{A}|\cdot \emph{\textsf{M}}_{\integ}(\log N))$ bit operations.}\label{c2}}

\begin{proof}
We modify the procedure of Algorithm \ref{a3}. After Steps 1-7, we additionally precompute $a^{x_1} \mod N$ and the values $a^{\Delta_{i,l_i}}\mod N$ for $i=1,\ldots,d$ and $l_i=1,\ldots,\kappa_i-1$. We may use a simple square-and-multiply algorithm for these modular exponentiations. Since all exponents are bounded by $M$, the runtime for the Steps 1-7 and these additional precomputations can be bounded by
$
O\left(\sum_{i=1}^d \kappa_i\cdot \textsf{M}_{\integ}(\log N)\log M\right)
$
bit operations, which is negligible due to our assumption. We now proceed with the main loop, with the only change that we replace the recursive formula for $x_{\nu+1}$ in Step 13 by
\[
a^{x_{\nu+1}}\mod N = a^{x_\nu}\cdot \prod_{i=1}^\mu a^{\Delta_{i,r_i}}\mod N.
\]
\end{proof}

\section{The sets $\mathcal{H}_{N,m}$ and $\mathcal{L}_{N,m}$\label{S2}}
In this section, we describe the core ideas for the proof of Theorem \ref{t3}. We will discuss results about the cardinality of certain sets containing the residues of divisors and the sums of divisors of natural numbers.

{\definition{Let $N,m\in\N$ such that $\gcd(N,m)=1$. We define the sets
\[
\mathcal{H}_{N,m}:=\{(x,y)\in\Z_m^2 :N\equiv xy \mod m\}
\]
and
\[
\mathcal{L}_{N,m}:=\{x+y\mod m: (x,y)\in\mathcal{H}_{N,m}\}.
\]}}

\remark{The set $\mathcal{H}_{N,m}$ is called a modular hyperbola. We refer the reader to \cite{Shp} for sophisticated results about the distribution of its elements. In this paper, we only use trivial properties:
\begin{enumerate}
\item{$|\mathcal{H}_{N,m}|=\phi(m)$.}
\item{If $N=uv$, then $(u\mod m, v \mod m)\in \mathcal{H}_{N,m}$.}
\end{enumerate}
In the following, we will focus on studying $\mathcal{L}_{N,m}$. Note that if $N=uv$, we clearly derive that $u+v\mod m\in\mathcal{L}_{N,m}$. For this reason, this set will be of interest in Section \ref{S1}.
}

{\lemma{
Let $N\in\N$ and $r$ be an odd prime such that $\gcd(N,r)=1$. Then
\[
  |\mathcal{L}_{N,r}|=
  \begin{cases}
    (r+1)/2      & \mbox{if } (N|r)=1,\\
    (r-1)/2      & \mbox{if } (N|r)=-1.
   \end{cases}
\]
Here, $(N|r)$ denotes the Legendre symbol.
}\label{l5}}

\begin{proof}
If $(N|r)=1$, then there are two distinct solutions to $x^2\equiv N \mod r$. Let $a$ and $b$ be these solutions. Then $(a,a)$ and $(b,b)$ are both elements of $\mathcal{H}_{N,r}$. Let $(x,y)\in\mathcal{H}_{N,r}$ be one of the other $r-3$ elements. Then we have $(y,x)\in\mathcal{H}_{N,r}$ and $(x,y)\neq (y,x)$. Since $x+y\equiv y+x \mod r$, we conclude that $(x,y)$ and $(y,x)$ pair up to one element in $\mathcal{L}_{N,r}$. Therefore, we deduce
\[
|\mathcal{L}_{N,r}|\leq \frac{r-3}{2}+2=\frac{r+1}{2}.
\]
If $(N|r)=-1$, we derive that $(x,y)\neq (y,x)$ holds for every element in $\mathcal{H}_{N,r}$. We obtain $|\mathcal{L}_{N,r}|\leq (r-1)/2$ by using similar arguments as above.

Now let $(x_0,y_0)$ and $(x_1,y_1)$ be two distinct elements of $\mathcal{H}_{N,r}$. Assuming that $x_0+y_0\equiv x_1+y_1 \mod r$ holds, we easily derive $N+y_0^2\equiv y_0(x_1+y_1) \mod r$ and, hence,
\[
y_0^2-y_0(x_1+y_1)+x_1y_1\equiv 0 \mod r.
\]
If we consider this as a quadratic equation in the variable $y_0$, the solutions are $x_1$ and $y_1$. Since the pairs are distinct, $y_0=y_1$ is not possible. Hence, we are left with $y_0=x_1$. We conclude that the only couples of elements of $\mathcal{H}_{N,r}$ which pair up to one element in $\mathcal{L}_{N,r}$ are the symmetrical ones we already considered. This finishes the proof.
\end{proof}

{\lemma{Let $N$ be a natural number and $m_1,\ldots,m_d\in\N$ be pairwise coprime. Furthermore, define $m:=m_1\cdots m_d$. Then $\mathcal{L}_{N,m}\subseteq \Z_m$ is the unique set determined by the $\mathcal{L}_{N,m_i}$ due to the CRT.}\label{l6}}

\begin{proof}
Consider the map
\begin{align*}
&\psi:\mathcal{L}_{N,m}\rightarrow \prod_{i=1}^{d} \mathcal{L}_{N,m_i},\\
x\mapsto (x&\mod m_1,\ldots,x\mod m_d).
\end{align*} 
We first show that $\psi$ is well-defined. Let $x\in\mathcal{L}_{N,m}$ and $i\in\{1,\ldots,d\}$ be arbitrarily chosen. By definition, there exists at least one pair $(x',y')\in\mathcal{H}_{N,m}$ such that $x\equiv x'+y' \mod m$. Clearly, $(x'\mod m_i, y'\mod m_i)$ is an element of $\mathcal{H}_{N,m_i}$. We deduce that the element $x \mod m_i=x'+y' \mod m_i$ is in $\mathcal{L}_{N,m_i}$, which we wanted to prove.

As a result of the CRT, $\psi$ is injective. We prove that this map is also surjective. Let $(s_1,\ldots,s_d)\in \prod_{i=1}^{d} \mathcal{L}_{N,m_i}$ be arbitrary. For every $i\in\{1,\ldots,d\}$, there exists a pair $(u_i,v_i)$ in $\mathcal{H}_{N,m_i}$ such that $s_i\equiv u_i+v_i\mod m_i$. By the CRT, we exhibit $u$ and $v$ such that
\[
u\equiv u_i \mod m_i \text{ and } v\equiv v_i \mod m_i
\]
and, hence, $N\equiv uv \mod m_i$ for every $i$. This implies that $N\equiv uv \mod m$, and $s:=u+v\mod m\in\mathcal{L}_{N,m}$. It is easy to see that $\psi(s)=(s_1,\ldots,s_d)$. Therefore, $\psi$ is bijective and the statement follows. 
\end{proof}

{\theorem{
Let $N\in\N$, $B\in\R^+$ and set $Q=\prod_{2<r\leq B} r$ for prime numbers $r$. Assume that $\gcd(N,Q)=1$. Then
\[
|\mathcal{L}_{N,Q}|\in O\left(\frac{Q\log B}{2^{\pi(B)}}\right),
\]
where $\pi(B)$ is the prime counting function. 
}\label{t7}}

\begin{proof}
As a result of Lemma \ref{l6}, we derive $|\mathcal{L}_{N,Q}|=\prod_{2<r\leq B}|\mathcal{L}_{N,r}|$. Note that Lemma \ref{l5} implies $|\mathcal{L}_{N,r}|\leq (r+1)/2$ for odd primes $r$. By using Corollary \ref{t5}, we conclude that
\[
\frac{|\mathcal{L}_{N,Q}|}{Q}\leq \frac{4}{3}\cdot\prod_{r\leq B}\frac{(r+1)/2}{r}\in O\left(\frac{\log B}{2^{\pi(B)}}\right).
\]
\end{proof}

\section{Order-finding algorithms}
In this section, we will discuss Sutherland's improvement of Shanks' babystep-giantstep method for the problem of finding the order of elements in $\Z_N^*$. Furthermore, we will present an algorithm that either finds a nontrivial factor of $N$, proves its primality or yields $a\in\Z_N^*$ whose order is larger than some given bound.

Let $E=1$ and $b=1$. Then Problem \ref{p1} is about finding the order of $a$ or proving that the order is larger than $T$. The following theorem shows that there is a way to asymptotically speed up Shanks' method by the factor $\sqrt{\log\log N}$.

{\theorem{There exists an algorithm with the following properties:
\begin{itemize}
\item{Input: $N\in\N$, $T\leq N$ and $a\in\Z_N^*$.}
\item{Output: If $\ord_N(a)\leq T$, then the output is $\ord_N(a)$; otherwise, the output is '$\ord_N(a)> T$'.}
\end{itemize}
The runtime complexity is bounded by
\[
O\Big{(}\frac{T^{1/2}}{\sqrt{\log\log T}}\cdot \log^2 N\Big{)}
\]
bit operations.}\label{t2}}

\begin{proof}
Assume that $\ord_N(a)\leq T$. We consider Algorithm 4.1. in \cite{Sut} for finding $\ord_N(a)$. The idea is to compute
\[
s:=a^{\prod_{r\leq B}r^{\lceil\log_r T\rceil}} \mod N
\]
for primes $r$ and a suitable choice of $B$. As a result, the order of $s$ is not divisible by any prime smaller than $B$. We may write $\ord_N(s)=mi+j$, where $m$ is of suitable size and divisible by 
\[
Q:=\prod_{r\leq B}r,
\]
and $j$ is coprime to $Q$. We proceed by computing $\ord_N(s)$. According to Proposition 4.2. in \cite{Sut}, we may compute suitable babysteps and giantsteps by performing $O(T^{1/2}(\log\log T)^{-1/2})$ multiplications modulo $N$. Combining Sutherland's result with the considerations in Remark \ref{r1}, we deduce the stated runtime for computing $\ord_N(s)$. Knowing the order of $s$, we may easily derive the order of $a$.

If this procedure does not find $\ord_N(a)$, it certainly follows that our assumption $\ord_N(a)\leq T$ is false.
\end{proof}

The following algorithm is a solution to the problem of either finding a nontrivial factor of a natural number, proving its primality or finding an element in its multiplicative group with order larger than some given bound. We will need such an element for applying Theorem \ref{t3} in order to prove Theorem \ref{A}. 

\begin{algorithm}
\emph{Input:} $N\in\N$ and $\delta\leq N$.

\emph{Output:} The algorithm either returns some $a\in\Z_N^*$ with $\ord_N(a)>\delta$, or some nontrivial factor of $N$, or '$N$ is prime'.

\begin{algorithmic}[1]
\State Set $M_1=1$ and $a=2$.
\For {$e=1,2,\ldots$}

	\While{$a\nmid N$ and $a^{M_e}\equiv 1 \mod N$}
		\State $a\gets a+1$
	\EndWhile
	\If{$a\mid N$}
		\State Return $a$ as nontrivial factor of $N$ or, if $a=N$, '$N$ is prime'.
	\EndIf
	\State Apply Theorem \ref{t2} with $T=\delta^{1/3}$.
	\If{$\ord_N(a)$ is not found}
		\State Apply Theorem \ref{t2} with $T=\delta$.
		\If{$\ord_N(a)$ is not found} 
			\State Return $a$ as element with $\ord_N(a)> \delta$.
		\EndIf
	\EndIf
	\State Set $m_{e}=\ord_N(a)$ and compute the prime factorization of $m_{e}$.
	\For{each prime $r$ dividing $m_{e}$}
		\If{$\gcd(N,a^{m_{e}/r}-1)\neq 1$}
			\State Return $\gcd(N,a^{m_{e}/r}-1)$ as nontrivial factor of $N$.
		\EndIf
	\EndFor
	\State $M_{e+1}\gets \lcm(M_e,m_{e})$
	\If{$M_{e+1}\geq \delta^{1/3}$}
		\State Apply Theorem \ref{t4} with $s=1$ and $m=M_{e+1}$.
		\State Return some nontrivial factor of $N$ or '$N$ is prime'.
	\EndIf
	\State $a\gets a+1$
\EndFor
\end{algorithmic}
\label{a1}
\end{algorithm}

{{\theorem{Algorithm \ref{a1} is correct. Assuming $N^{2/5}\leq\delta$, its runtime complexity is bounded by
\[
O\left(\frac{\delta^{1/2}}{\sqrt{\log\log \delta}}\cdot \log^2 N \right)
\]
bit operations.}\label{t1}}

\begin{proof}
\textit{Correctness:}
We first show that, if the algorithm terminates, it yields a desired output. If the algorithm stops in Step 6, 11 or 15, it is easy to see that we have either determined primality, a nontrivial factor of $N$ or found $a\in\Z_N^*$ such that $\ord_N(a)> \delta$. Let us assume that the algorithm has reached Step 16. Since no proper factor has been found in the loop from Step 13 to 15, Lemma \ref{l2} yields that $m_{e}=\ord_p(a)$ for every prime factor $p$ of $N$. Note that this implies $p\equiv 1 \mod m_{e}$. Furthermore, this congruence also holds modulo $m_j$ for $j=1,\ldots,e-1$, because if otherwise, $N$ would have been factored in one of the previous runs of the loop. We conclude that the congruence $p\equiv 1 \mod M_{e+1}$ holds for every prime factor $p$ of $N$. Hence, Theorem \ref{t4} is correctly applied in Step 18.

We are left with proving that the algorithm eventually terminates. First note that at the entry of the main loop in Step 2, it certainly holds that $M_e$ is a divisor of $\phi(N)$, that $M_e<\delta^{1/3}$ and that $b^{M_e}\equiv 1\mod N$ for each $b<a$. We proceed by showing that the while loop in Step 3 and 4 terminates after at most $M_e$ runs. Assume to the contrary that this loop runs for $M_e+1$ different values of $a$. The smallest prime factor $p$ of $N$ is always larger than $a$, since $N$ has not been factored in previous runs. Furthermore, these elements all satisfy $a^{M_e}\equiv 1  \mod N$, which implies  $a^{M_e}\equiv 1 \mod p$. But since $p$ is prime, there are at most $M_e$ elements in $\Z_p^*$ satisfying this property. As a consequence, we derive a contradiction.

We now want to prove that the main loop terminates after at most $O(\log \delta)$ runs. Consider the $e$-th run of the loop. If no factor is found in Step 5, we have $a\nmid N$ and $a^{M_{e}}\not\equiv 1 \mod N$ in Step 7. If the algorithm does not terminate in Step 11, we have found the order of $a$ in Step 7 or Step 9. Note that $m_{e}$ does not divide $M_{e}$ and, hence, the value of $M_{e+1}$ computed in Step 16 is at least twice as large as $M_{e}$. As a consequence, this value exceeds $\delta^{1/3}$ after at most $O(\log \delta)$ runs and, hence, the algorithm eventually terminates.

\textit{Running time:}
We start by considering the running time of Step 18. According to Theorem \ref{t4}, the cost for this step can be bounded by
\[
O\left(\left(\frac{N^{1/4}}{\sqrt{M_e}}\right)^{1+o(1)}\right)
\]
bit operations. Since the algorithm reached Step 18, it follows that $M_e\geq \delta^{1/3}$. Furthermore, we assumed $N^{2/5}\leq \delta$. We conclude
\[
N^{1/4}/\sqrt{M_e} \leq N^{1/4}/\delta^{1/6}\leq\delta^{11/24}.
\]
Hence, the cost for performing Step $18$ is asymptotically negligible.

We continue by showing that, once the algorithm reaches Step 9, it terminates within the claimed running time. Using Theorem \ref{t2}, Step 9 can be performed in 
\[
O\left(\frac{\delta^{1/2}}{\sqrt{\log\log \delta}}\cdot \log^2 N\right)
\]
bit operations. If we do not find the order of $a$, the algorithm terminates. We assume that $\ord_N(a)\leq \delta$ is found. In Step 12, we may compute the divisors of $m_e$ in negligible time by using any factorization algorithm which is at least slightly faster than naive trial division. The cost for computing the greatest common divisors in Step 14 and the value of $M_e$ in Step 16 is negligible. Now note that since the algorithm reached Step 9, it follows that $m_e> \delta^{1/3}$. Therefore, we derive $M_e\geq m_e > \delta ^{1/3}$ and the algorithm reaches Step 18, whose running time has already been analyzed and is asymptotically negligible.

We now consider the complexity of the $e$-th run of the main loop. One run of the while loop in the Steps 3 and 4 can be performed in polynomial time and we have already shown that there are at most $M_{e}<\delta^{1/3}$ runs of this loop. Therefore, its overall running time is bounded by $O(\delta^{1/3+o(1)})$ bit operations. By Theorem \ref{t2}, we may perform Step 7 in $O(\delta^{1/6+o(1)})$ bit operations. If $\ord_N(a)$ is not found, the algorithm reaches Step 9 and is going to terminate within the claimed running time. We assume that $\ord_N(a)$ is found. Since $m_{e}<\delta ^{1/3}$, we may use trial division in Step 12 to compute the divisors of $m_{e}$ in $O(\delta^{1/6+o(1)})$ bit operations. The computations in the Steps 14 and 16 can be performed in polynomial time. 

We derive that the overall bit-complexity of one complete run of the main loop is bounded by $O(\delta^{1/3+o(1)})$. Since we have already shown that there are at most $O(\log \delta)$ runs of this loop, we conclude that the algorithm reaches either Step 9 or Step 18 in asymptotically negligible time. This finishes the proof.
\end{proof}

\remark{Concerning Algorithm \ref{a1}, we would like to point out the following: 

\begin{enumerate}
\item{In order to prove Theorem \ref{A}, it suffices to suppose $N^{2/5}\leq \delta$. However, the arguments in the proof of the running time of the algorithm have shown that this assumption is not optimized. By varying the exponent of $\delta^{1/3}$ in the Steps 7 and 17 and by comparing the logarithmic factors of the complexity bounds of Theorem \ref{t4}, Theorem \ref{t2} and the main loop, one should be able to prove the result with an assumption of the form $\delta\in\Omega(N^{1/3+o(1)})$.}
\item{Note that the purpose of Algorithm \ref{a1} is entirely theoretical. The procedure will not find a factor of $N$ in most cases, at least not for reasonable values of $\delta$. It should be considered as a deterministic and rigorous method for solving the underlying computational problem and finding an element of sufficiently large order.}
\end{enumerate}}

\section{On solving $a^X\equiv a^{N+1} \mod N$\label{S1}}
The following theorem yields a way to improve the standard BSGS technique for solving the discrete logarithm problem mentioned in the title. This result will be substantial for proving Theorem \ref{A} in Section \ref{S4}.

{\theorem{There exists an algorithm with the following properties: 
\begin{itemize}
\item{Input: $N\in\N$, $N^{1/2}\leq T\leq N$ and $a\in\Z_N^*$ such that $\ord_N(a)\geq (T/2)^{2/3}$.}
\item{Output: A set $\mathcal{L}$ of integers in $[1,O(T)]$ with $|\mathcal{L}|\in O(T^{1/2-\xi(T)}\sqrt{\log\log T})$, where
\[
\xi(T):=\frac{\log \sqrt[4]{2}}{\log(\frac{1}{2}\log T)}.
\]}
\end{itemize}
Its bit-complexity is bounded by
\[
O(T^{1/2-\xi(T)}\sqrt{\log\log T}\cdot \log^2 N).
\]
In the case when $N$ is semiprime with prime factors $p$ and $q$ and $p+q\leq T$, we have $p+q\in\mathcal{L}$. 
}\label{t3}}

\begin{proof}
Assume that $N$ is semiprime with prime factors $p$ and $q$ and suppose that $p+q=:S\leq T$. We compute $b:=a^{N+1}\mod N$, which can be done in polynomial time. By Lemma \ref{l3} we have $a^S\equiv b \mod N$. Therefore, we want to find a set $\mathcal{L}$ of solutions to
\[
a^x\equiv b \mod N
\]
such that $S\in\mathcal{L}$. For this task, we use the results of Section \ref{S2} on the residue classes containing $S$ to adapt Shanks' babystep-giantstep method. Discussing asymptotic complexities, we will assume that the input $T$ is sufficiently large. Let $B=\frac{1}{2}\log T$. For prime numbers $r$ and Chebyshev's theta function $\vartheta(x)$, define
\[
Q:=\prod_{2<r\leq B}r=\frac{1}{2}\cdot \prod_{r\leq B}r=\frac{e^{\vartheta(B)}}{2}.
\]
We first want to show that $Q$ is smaller than $T^{1/2}\cdot(\log B)^{-1/2}\cdot 2^{\pi(B)/2}$. Considering the expression above, we easily deduce that it suffices to prove that the inequality $e^{\vartheta(B)}<T^{1/2}\cdot(\log B)^{-1/2}\cdot 2^{\pi(B)/2}$ holds, which is equivalent to
\[
\vartheta(B) <B-\frac{1}{2}\log \log B+\frac{\log 2}{2}\pi(B).
\]
Note that this inequality is true due to Lemma \ref{l7}, which proves the claim. Next, we compute a numerical approximation $\alpha$ such that $|T^{1/2}\cdot(\log B)^{-1/2}\cdot 2^{\pi(B)/2}-\alpha|<1$. Find $k\in\N$ such that $2^{k-1}Q<\alpha +1\leq 2^{k}Q$ and define $m:=2^kQ$. Obviously, we may derive the value of $m$ in polynomial time. Now perform trial division up to $B$ to check if one of the prime divisors of $N$ is smaller than $B$. If this is the case, we have found $S$ itself and the problem is solved. Assume that trial division has failed. Then we know that $N$ is coprime to $m$. Since $S\leq T$, we have $S=mu+v$ for $v:=S\mod m$ and some $u$ with
\[
u\in\{0,1,\ldots,\lceil T^{1/2}\cdot\sqrt{\log B}\cdot 2^{-\pi(B)/2}\rceil-1\}.
\]
Let $\beta $ be such that $|T^{1/2}\cdot\sqrt{\log B}\cdot 2^{-\pi(B)/2}-\beta|<1$. For $j=0,1,\ldots,\lceil \beta \rceil$, we compute the giantsteps of the form $a^{-mj}b\mod N$, which can be done in
\begin{equation}
O\left(\frac{T^{1/2}\cdot\sqrt{\log B}}{2^{\pi(B)/2}}\cdot\textsf{M}_{\integ}(\log N)\right)=O\Big{(}\frac{T^{1/2}\sqrt{\log\log T}}{2^{\pi(B)/2}}\cdot \textsf{M}_{\integ}(\log N)\Big{)}
\end{equation}
bit operations by precomputing the value $a^{-m}\mod N$ and using the recursive formula $a^{-m(j+1)}b\equiv a^{-m}a^{-mj}b \mod N$.

Next, we discuss the computation of the babysteps, where one of them should be equal to the value $a^v\mod N$. We proceed as follows: Note that Lemma \ref{l6} implies that $\mathcal{L}_{N,m}$ is determined by $\mathcal{L}_{N,2^k}$ and $\mathcal{L}_{N,r}$, $2<r\leq B$. We start by computing these sets. Let $r$ be an arbitrary prime with $2<r\leq B$. In order to find the elements of the set $\mathcal{L}_{N,r}$, we first determine $\mathcal{H}_{N,r}$. This can be done easily by evaluating $(x,x^{-1}N\mod r)$, where $x=1,2,\ldots,r-1$. We then find the elements of $\mathcal{L}_{N,r}$ via the definition of these sets. Furthermore, the elements of $\mathcal{L}_{N,2^k}$ may be computed in a similar manner. For applying Corollary \ref{c2}, we have to ensure that these sets are stored on tape with decreasing cardinalities. We first compute $\mathcal{L}_{N,2^k}$ and $|\mathcal{L}_{N,2^k}|$. By Lemma \ref{l5}, we deduce that $r_1>r_2$ implies $|\mathcal{L}_{N,r_1}|\geq|\mathcal{L}_{N,r_2}|$ for odd primes $r_1,r_2$. Hence, we proceed by computing $\mathcal{L}_{N,r}$ and $|\mathcal{L}_{N,r}|$ for $2<r\leq B$, starting with the set corresponding to the largest prime in this range. By comparing $|\mathcal{L}_{N,2^k}|$ to $(r-1)/2$ for $2<r\leq B$, we may put the elements of $\mathcal{L}_{N,2^k}$ in the right place. The total runtime for these tasks is bounded by $O(2^k+\sum_{r\leq B} r)$ arithmetic operations. Due to our assumption of $T$ being sufficiently large, the cost for these computations is negligible and the assumption of Corollary \ref{c2} is satisfied. Hence, we may use this result to compute the elements $a^i\mod N$ for $i\in\mathcal{L}_{N,m}$ in 
$
O(|\mathcal{L}_{N,m}|\cdot\textsf{M}_{\integ}(\log N))
$
bit operations. Since we have $m=2^kQ$ and since $2^k$ and $Q$ are coprime, by Lemma \ref{l6} it also follows that $|\mathcal{L}_{N,m}|=|\mathcal{L}_{N,2^k}|\cdot|\mathcal{L}_{N,Q}|$. Note that $|\mathcal{L}_{N,2^k}|\in O(2^k)$. By Theorem \ref{t7}, we derive that
\[
|\mathcal{L}_{N,m}|\in O\left(\frac{m\log B}{2^{\pi(B)}}\right)\subseteq O(T^{1/2}\cdot \sqrt{\log B}\cdot 2^{-\pi(B)/2}).
\]
Since $v\in\mathcal{L}_{N,m}$, we conclude that we have computed suitable babysteps and the running time for doing so was bounded by (7.1). 

Set $\mathcal{L}=\emptyset$. We now apply merge sort to the lists of babysteps and giantsteps. As we have suggested in Remark \ref{r1}, we then search for matches. If the giantstep $a^{-mj}b\mod N$ is equal to the babystep $a^i\mod N$, we add the value $mj+i$ to $\mathcal{L}$. Since we have assumed $\ord_N(a)\geq (T/2)^{2/3}$ and since it is easy to check that $(T/2)^{2/3}>m$ holds for sufficiently large $T$, the babysteps are all distinct. As a consequence, not more than one value is added to $\mathcal{L}$ for every $j$, which yields
\[
|\mathcal{L}|\in O(T^{1/2}\cdot 2^{-\pi(B)/2}\cdot \sqrt{\log\log T}).
\]
In view of the discussion in Remark \ref{r1}, we conclude that the number of bit operations to finish this procedure is bounded by
\[
O\left(\frac{T^{1/2}\cdot\sqrt{\log\log T}}{2^{\pi(B)/2}}\log^2 N\right).
\]
Note that $\mathcal{L}$ contains $S$ as an element. We are left with the task to discuss the term $2^{-\pi(B)/2}$ in the bound. By the Prime Number Theorem, we have $\pi(B) \sim B/\log B$. The first statement of Theorem \ref{t8} implies that we have $B/\log B-2<\pi(B)$ for $B\geq 2$. We conclude that
\[
O\left(2^{-\frac{\pi(B)}{2}}\right)\subseteq O\left(2^{-\frac{B}{2\log B}}\right)=O\left(2^{-\frac{\frac{1}{4}\log T}{\log\left(\frac{1}{2}\log T\right)}}\right)=O\left(T^{-\xi(T)}\right).
\]
This finishes the proof.
\end{proof}

\section{Proof of Theorem \ref{A}\label{S4}}
Suppose that we want to compute the prime factorization of $N\in\N$. First, we perform the following preparation step: We want to use Corollary \ref{c3} and have to ensure that all divisors $d$ of $N$ satisfy $4d^{1/3}\leq d$. Hence, we use trial division up to $7$ to remove small prime factors from $N$. Next, we apply Corollary \ref{c3} to the resulting number $N_0$, putting $\delta_0=\lceil N_0^{1/3}\rceil$. If any prime divisor of $N_0$ is found, we remove it and denote the resulting number by $N_1$. We then apply Corollary \ref{c3} to $N_1$, setting $\delta_1=\lceil N_1^{1/3}\rceil$. We proceed in this way until no more divisors are found. Note that the number of prime factors of $N$ is bounded by $O(\log N)$. Since $\delta_i\leq \lceil N^{1/3}\rceil$ for every $i\in\N_0$, the bit-complexity of this procedure is bounded by
\[
O\left(\textsf{M}_{\integ}\Big{(}\frac{N^{1/6}\log N}{\sqrt{\log\log N}}\Big{)}\log N\right).
\]
Compared to the complexity bound stated in Theorem \ref{A}, this is asymptotically negligible. Let $N_k$ be the number which remains after these computations. Since no factor smaller or equal to $\lceil N_k^{1/3} \rceil$ has been found by applying Corollary \ref{c3}, $N_k$ does not have more than two prime factors. We may check easily if $N_k$ is a square number; if not, $N_k$ has to be a prime or a semiprime number. 

As a result of this argument, we may restrict our attention to the case where $N$ is a prime or a semiprime number. We now consider the following algorithm.
\begin{algorithm}
\emph{Input:} A prime or semiprime $N\in\N$ and a parameter $\Delta\in\R$ satisfying $N^{2/5}\leq\Delta\leq N^{1/2}$.

\emph{Output:} The prime factorization of $N$.

\begin{algorithmic}[1]
\State Apply Corollary \ref{c3} with $\delta=\Delta$. \Comment{Use Strassen's algorithm} 

If a prime factor of $N$ is found, return and stop. 
\State Apply Theorem \ref{t1} with $\delta=\Delta$.  

If a prime factor of $N$ is found or $N$ is proven to be prime, return and stop; otherwise, $a\in\Z_N^*$ with $\ord_N(a)>\Delta$ is found.
\State Apply Theorem \ref{t3} with $T=N^{1/2}+N/\Delta$. \Comment{Solve $a^X\equiv a^{N+1}\mod N$}

 Let $\mathcal{L}$ be the set computed by the algorithm.
\State For every $s\in \mathcal{L}$, apply the procedure of Lemma \ref{l1} to $N$ and $s$.

If a prime factor of $N$ is found, return and stop; otherwise, return that $N$ is prime.
\end{algorithmic}
\label{a2}
\end{algorithm}

For analyzing the runtime complexity of Algorithm \ref{a2}, we consider the notion of logarithmic multiplication time functions $\textsf{e}_{\integ}(k)$ defined in the introduction. We already know that $\log k\cdot 8^{\log^*k}$ is the asymptotically smallest known example for such a function, where $\log^*k$ denotes the iterated logarithm:
\[
  \log^* k :=
  \begin{cases}
    0                  & \mbox{if } k \leq 1, \\
    1 + \log^*(\log k) & \mbox{if } k > 1.
   \end{cases}
\]
Harvey, van der Hoeven and Lecerf deduced this result in \cite{HarHoeLec} by improving F\"urer's multiplication algorithm \cite{Für}.

{\lemma{Algorithm \ref{a2} is correct. Its runtime complexity is bounded by
\[
\max\Big{\{}O\Big{(}\emph{\textsf{M}}_{\integ}\Big{(}\frac{\sqrt{\Delta}\log N}{\sqrt{\log\log N}}\Big{)}\Big{)}, O\Big{(}\left(\frac{N}{\Delta}\right)^{\frac{1}{2}-\xi\left(\frac{N}{\Delta}\right)}\sqrt{\log\log N}\cdot\log^2 N\Big{)}\Big{\}}
\]
bit operations, where 
\[
\xi(x):=\frac{\log \sqrt[4]{2}}{\log(\frac{1}{2}\log x)}.
\]
}\label{l4}}

\begin{proof}
\textit{Correctness}: Assume that $N$ is semiprime and has the prime factors $p$ and $q$, where $p<q$. If $p$ is not found in Step 1, then $p>\Delta$. We deduce $q<N/\Delta$, which implies $S:=p+q<N^{1/2}+N/\Delta$. Since we have assumed $N^{2/5}\leq \Delta$, we are able to  apply Theorem \ref{t1} in Step 2. Assume that $N$ does not get factored. In Step 3, we have 
\[
T=N^{1/2}+N/\Delta\leq N^{1/2}+N^{3/5}\leq 2 N^{3/5},
\]
and we obtain
\[
\ord_N(a)>\Delta\geq N^{2/5}\geq (T/2)^{2/3}.
\]
Since we have already seen that $S\leq T$ holds, the set $\mathcal{L}$ computed in Step 3 contains $S$ as an element. Therefore, $N$ does get factored in Step 4. 

We have proven that the algorithm will determine a factor of $N$ if the number is semiprime. Therefore, if this procedure fails to find a factor, $N$ must be prime. We conclude that the algorithm is correct.

\textit{Running time}: The running time of Step 1 is in 
\[
O\Big{(}\textsf{M}_{\integ}\Big{(}\frac{\sqrt{\Delta}\log N}{\sqrt{\log\log N}}\Big{)}\Big{)}
\]
bit operations. Note that $\Delta \geq N^{2/5}$ yields $O((\log\log \Delta)^{-1/2})\subseteq O((\log\log N)^{-1/2})$.  Hence, the bit-complexity of Step 2 may also be bounded by
\[
O\Big{(}\frac{\Delta^{1/2}}{\sqrt{\log\log \Delta}}\cdot \log^2 N\Big{)}\subseteq O\Big{(}\textsf{M}_{\integ}\Big{(}\frac{\sqrt{\Delta}\log N}{\sqrt{\log\log N}}\Big{)}\Big{)}.
\] 
Here, we have used (2) and (3) of Definition \ref{d1}. Now since $N^{2/5}\leq \Delta \leq N^{1/2}$, it is easy to see that the bit-complexity of Step 3 is bounded by
\[
O\left(\left(\frac{N}{\Delta}\right)^{\frac{1}{2}-\xi\left(\frac{N}{\Delta}\right)}\sqrt{\log\log N}\cdot\log^2 N\right).
\]
We consider the running time of Step 4. Theorem \ref{t3} also yields that the set $\mathcal{L}$ computed in Step 3 consists of at most 
\[
O\left(\left(\frac{N}{\Delta}\right)^{\frac{1}{2}-\xi\left(\frac{N}{\Delta}\right)}\sqrt{\log\log N}\right)
\]
elements. Since the cost for one run of Lemma \ref{l1} is in $O(\textsf{M}_{\integ}(\log N))$ bit operations, this concludes the proof.
\end{proof}

To finish the proof of Theorem \ref{A}, we want to find an optimal choice for the value $\Delta$ such that the bound shown in Lemma \ref{l4} is minimal. First note that we have supposed that $N^{2/5} \leq \Delta \leq N^{1/2}$. The choice $\Delta=N^{1/2}$ yields the bound shown by Costa and Harvey. Now assume $\Delta<N^{1/2-c}$ for some constant $c$. Then we derive $N/\Delta>N^{1/2+c}$, and since $\xi(N/\Delta)$ tends to $0$ for increasing values of $N$, it is easy to see that the resulting value of the maximum in Lemma \ref{l4} is asymptotically worse compared to the bound of Costa and Harvey. As a result, we know that the optimal choice is of the form $\Delta=N^{1/2-\varepsilon(N)}$ for some $\varepsilon(N)\in o(1)$. 

Let $\varepsilon>0$ be arbitrary and define $\varepsilon':=2^{5/4}\varepsilon$. Note that our last observation implies that we may assume $\Delta\geq N^{1/2-\varepsilon'}$ for sufficiently large input $N\geq C_\varepsilon$. Our goal now is to determine an explicit value for $\Delta$. In order to do this, we want to solve
\begin{equation}
\frac{\sqrt{\Delta}\log N}{\sqrt{\log\log N}}\cdot \textsf{e}_{\integ}(N)=\left(\frac{N}{\Delta}\right)^{\frac{1}{2}-\xi\left(N^{1/2+\varepsilon'}\right)}\sqrt{\log\log N}\cdot \log^2 N
\end{equation}
for $\Delta$. Ignoring implied constants, the left-hand side of $(8.1)$ is the first argument of the maximum in Lemma \ref{l4}, where we have applied (1) and (2) of Definition \ref{d1}. The right-hand side of $(8.1)$ is an \emph{upper bound} of the second argument of this maximum, where we have replaced $\xi(N/\Delta)$ by $\xi(N^{1/2+\varepsilon'})$. Comparing these two sides of the equation, the resulting value of $\Delta$ is optimal; because if we vary $\Delta$, the left-hand side is increasing if and only if the right-hand side is decreasing. Equation $(8.1)$ can be rewritten as
\[
\Delta=N^{\frac{1/2-\xi\left(N^{1/2+\varepsilon'}\right)}{1-\xi\left(N^{1/2+\varepsilon'}\right)}}\cdot (R_{N})^{\frac{1}{1-\xi\left(N^{1/2+\varepsilon'}\right)}},
\]
where
\[
R_N:=\frac{\log\log N\cdot  \log N}{\textsf{e}_{\integ}(N)}.
\]
We consider these logarithmic factors first. Note that the behaviour of $R_N$ depends on $\textsf{e}_{\integ}$. By using F\"urers bound for this function, we immediately derive $R_N\in\Omega(1)$ and the rough upper bound $R_N<\log^2 N$. In any case, we obtain
\[
R_N^{\frac{1}{1-\xi\left(N^{1/2+\varepsilon'}\right)}}=R_N \cdot R_N ^{\frac{\xi\left(N^{1/2+\varepsilon'}\right)}{1-\xi\left(N^{1/2+\varepsilon'}\right)}}.
\]
Since $\xi\left(N^{1/2+\varepsilon'}\right)$ tends to $0$ for $N$ to infinity, $R_N\in\Omega(1)$ clearly implies
\begin{align*}
R_N ^{\frac{\xi\left(N^{1/2+\varepsilon'}\right)}{1-\xi\left(N^{1/2+\varepsilon'}\right)}}\in\Omega(1).
\end{align*}
Furthermore, one easily deduces
\[
R_N ^{\frac{\xi\left(N^{1/2+\varepsilon'}\right)}{1-\xi\left(N^{1/2+\varepsilon'}\right)}}<(\log N)^{\frac{2\xi\left(N^{1/2+\varepsilon'}\right)}{1-\xi\left(N^{1/2+\varepsilon'}\right)}}=(\log N)^\frac{\log \sqrt{2}}{\log\left(2^{-5/4}\log \left(N^{1/2+\varepsilon'}\right)\right)}\in O(1).
\]
As a consequence of these estimates, we conclude 
\[
R_N^{\frac{1}{1-\xi\left(N^{1/2+\varepsilon'}\right)}}\in\Theta (R_N).
\]
We now consider the term
\[
N^{\frac{1/2-\xi\left(N^{1/2+\varepsilon'}\right)}{1-\xi\left(N^{1/2+\varepsilon'}\right)}}=N^{\frac{1}{2}-\frac{\xi\left(N^{1/2+\varepsilon'}\right)}{2-2\xi\left(N^{1/2+\varepsilon'}\right)}}
\]
and derive
\[
\frac{\xi(N^{1/2+\varepsilon'})}{2-2\xi(N^{1/2+\varepsilon'})}=\frac{\log \sqrt[8]{2}}{\log\left(2^{-5/4}\log (N^{1/2+\varepsilon'})\right)}=\frac{\log \sqrt[8]{2}}{\log\left((2^{-9/4}+\varepsilon)\log N\right)}.
\]
Hence, any
\[
\Delta\in\Theta\left(N^{\frac{1}{2}-\frac{\log{2}}{8\log\left((2^{-9/4}+\varepsilon)\log N\right)}} \cdot R_N\right)
\]
is an asymptotically optimal value. Since $\varepsilon$ was arbitrary, substituting such $\Delta$ on the left-hand side of $(8.1)$ finally yields the following result.

{\theorem{Let $\varepsilon >0$ be arbitrary and $N$ be any natural number. Then the prime factorization of $N$ can be computed in
\[
O\left(N^{\frac{1}{4}}\exp\left(\frac{-C\log N}{\log\left((2^{-9/4}+\varepsilon)\log N\right)}\right)\cdot\sqrt{\emph{\textsf{e}}_{\integ}(N)\log^3 N}\right)
\]
bit operations, where $C=(\log 2)/16$ and the only implied constant depends on $\varepsilon$.}\label{t9}}

\remark{For any $\varepsilon>0$, we may compute an optimal value $\Delta_\varepsilon$ via the formula given above. We have already seen that our assumption $\Delta_\varepsilon\geq N^{1/2-\varepsilon}$ is satisfied if $N$ is larger than some computable bound $C_\varepsilon$. The bound $C_\varepsilon$ is increasing if $\varepsilon$ tends to $0$.

Setting $\varepsilon=0.2103-2^{-9/4}\approx 7.59\cdot 10^{-5}$ in Theorem \ref{t9} yields the statement of Theorem \ref{A} and finishes the proof.\label{r4}}

\section{Addendum}
We address two problems concerning refinements of runtime complexities.
\problem{In the proof of Theorem \ref{t3}, the choice of $B=\frac{1}{2}\log T$ affected the runtime complexity bound. Therefore, the question arises whether this choice is optimal. The advantage of our method mainly depends on the increase of $\pi(B)$. Hence, $B$ should be chosen as large as possible. Following the arguments in the proof, one observes that we cannot set $B=\frac{1}{c}\log T$ for any $c<2$. In this sense, our choice of $B$ is optimal. However, there may be some $c(T)$ converging to $0$ such that $B=(2-c(T))^{-1}\log T$
works out for sufficiently large $T$. At this point, it remains an open question whether we may improve the complexity bound with this idea.}

\problem{We consider the analysis of the running time of Algorithm \ref{a2} in Section \ref{S4}. For any choice of $\varepsilon$, we use an upper bound and not the actual runtime complexity on the right-hand side of equation $(8.1)$ to determine an explicit value for $\Delta$. Therefore, the resulting exponent in the runtime complexity is close to optimal for small $\varepsilon$, but never best possible. We have started with the observation that $\Delta=N^{1/2-\varepsilon(N)}$ holds for some $\varepsilon(N)\in o(1)$. Replacing $\varepsilon$ by $\varepsilon(N)$ in our estimates, the resulting exponent is of the form
\[
\frac{1}{2}-\frac{\log \sqrt[8]{2}}{\log\left((2^{-9/4}+\varepsilon(N))\log N\right)}.
\]
Instead of starting with $\varepsilon$ and computing the lower bound $C_\varepsilon$ for $N$ discussed in Remark \ref{r4}, one could go the other way and ask for given $N$, what the corresponding value $\varepsilon=\varepsilon(N)$ is. Unfortunately, determining $\varepsilon(N)$ explicitly and deducing a closed formula for the general optimal value of $\Delta$ seems to be infeasible.}
\vspace{12pt}
\parindent 0pt

\textbf{Acknowledgement.} Special thanks go to an anonymous referee for an extensive and most valuable report, which contributed to improve the exposition of this paper.



\begin{thebibliography}{HD}
\bibitem[BGS07]{BosGauSch} A. Bostan, P. Gaudry, \'{E}. Schost, \emph{Linear recurrences with polynomial coefficients and application to integer factorization and Cartier-Manin operator}, SIAM J. Comput., 36(6): 1777-1806, 2007.

\bibitem[CH14]{CosHar} E. Costa, D. Harvey,
\emph{Faster deterministic integer factorization}, Math. Comp., 83: 339-345, 2014.

\bibitem[F\"ur09]{Für} M. F\"urer,
\emph{Faster integer multiplication}, SIAM J. Comput., 39(3): 979-1005, 2009.

\bibitem[HHL16]{HarHoeLec} D. Harvey, J. v. d. Hoeven, G. Lecerf,
\emph{Even faster integer multiplication}, J. Complexity, 36: 1-30, 2016.

\bibitem[Hit15]{Hit} M. Hittmeir,
\emph{Deterministic factorization of sums and differences of powers}, Math. Comp., to appear, 2015.

\bibitem[Len00]{Len} A. K. Lenstra,
\emph{Integer Factoring}, Designs, Codes and Cryptography, 19: 101-128, 2000.

\bibitem[Rie94]{Rie} H. Riesel,
\emph{Prime Numbers and Computer Methods for Factorization}, Progress in Mathematics (Volume 126), Second Edition, Birkh\"auser Boston, 1994.

\bibitem[RS62]{RosSch} J.B. Rosser, L. Schoenfeld,
\emph{Approximate formulas for some functions of prime numbers}, Illinois J. Math., 6: 64-94, 1964.

\bibitem[Sch76]{Sch} L. Schoenfeld,
\emph{Sharper Bounds for the Chebyshev Functions $\theta(x)$ and $\psi(x)$. II}, Math. Comp., 30: 337-360, 1976.

\bibitem[SW11]{SedWay} R. Sedgewick, K. Wayne,
\emph{Algorithms}, Princeton University, Fourth Edition, Addison-Wesley, 2011.

\bibitem[Sho05]{Sho} V. Shoup,
\emph{A Computational Introduction to Number Theory and Algebra}, Cambridge University Press, 2005.

\bibitem[Shp12]{Shp} I. E. Shparlinski,
\emph{Modular hyperbolas}, Jap. J. Math., 7: 235-294, 2012.

\bibitem[Str77]{Str} V. Strassen,
\emph{Einige Resultate \"uber Berechnungskomplexit\"at},
Jber. Deutsch. Math.-Verein., 78(1): 1-8, 1976/77.

\bibitem[Sut07]{Sut} A. V. Sutherland,
\emph{Order Computations in Generic Groups},
Doctoral Thesis, Massachusetts Institute of Technology, 2007.

\bibitem[Wag13]{Wag} S.S. Wagstaff Jr.,
\emph{The Joy of Factoring},  
American Math. Society, Providence, RI, 2013.
\end{thebibliography}
\end{document}